\documentstyle[12pt]{amsart}

\newcommand{\bth}[1]{\bigskip {\bf #1} \it}
\newcommand{\eth}{\rm \medskip}

\newfont{\bbb}{msbm10}

\newcommand{\R}{\mbox{\bbb R}}

\newtheorem{lem}{Lemma}

\newtheorem{prop}{Proposition}

\title{Curvature properties of the Chern connection of twistor spaces}
\author{ Johann Davidov, Gueo Grantcharov and Oleg  Mu\u skarov}

\address{Johann Davidov\\ Institute of Mathematics and Informatics \\
Bulgarian Academy of Sciences\\ Acad. G.Bonchev Str. Bl.8\\ 1113 Sofia\\ Bulgaria} \email{jtd@@math.bas.bg}
\address{Gueo Grantcharov\\ Department of Mathematics\\
Florida International University\\ Miami\\FL 33199\\USA}
\email{grantchg@@fiu.edu}
\address{Oleg Mu\u skarov \\Institute of Mathematics and Informatics \\
Bulgarian Academy of Sciences\\ Acad. G.Bonchev Str. Bl.8\\ 1113 Sofia\\ Bulgaria}
\email{muskarov@@math.bas.bg}
\thanks{The authors are partially supported by  EDGE Research Training Network HPRN-CT-2000-00101.
 The second author is supported also by  NSF Grant DMS 0333172}

\date{}

\begin{document}

\maketitle

\begin{abstract}
  The twistor space ${\cal Z}$ of an oriented Riemannian $4$-manifold $M$ admits a natural
  $1$-parameter family of
 Riemannian metrics $h_t$ compatible with the almost complex structures $J_1$ and $J_2$ introduced, respectively,
by Atiyah, Hitchin and Singer, and Eells and Salamon. In this paper we compute the first Chern form of the almost
Hermitian manifold $({\cal Z},h_{t},J_{n})$, $n=1,2$ and find the geometric conditions on $M$ under which the
curvature of its Chern connection $D^{n}$ is of type $(1,1)$. We also describe the twistor spaces of constant
holo%
morphic sectional curvature with respect to $D^n$ and show that the Nijenhuis tensor of $J_2$ is $D^2$-parallel
provided the base manifold $M$ is Einstein and self-dual.
\end{abstract}

\vspace{0,5cm} \hspace{0,7cm}{\it Keywords}: Twistor spaces, Chern connection

\vspace{0,2cm} \hspace{0,7cm}{\it 1991 MSC: 53C15}

\parindent 1cm

\section{Introduction}

     It is well-known \cite{Lich, Gaud} that every almost Hermitian manifold admits a unique connection
for which the almost complex structure and the metric are parallel and the $(1,1)$-part of the torsion vanishes.
It is usually called the Chern connection because, in the integrable case, it coincides with the Chern connection
\cite{Chern} of the tangent bundle considered as a Hermitian holomorphic bundle. This connection plays an
important role in (almost) complex geometry since, by the Chern- Weil theory, the Chern classes of the manifold
are directly  related to its curvature.

    Motivated by the recent works of S.Donaldson \cite{Don} and C.LeBrun \cite{LB}, V.Apos-\\tolov and T.Dragichi
    \cite{AD}
proposed to study the problem of existence of almost-K\"ahler
structures of constant Hermitian scalar curvature and/or type
$(1,1)$ Ricci form of its Chern connection (from now on we refer
to it as the first Chern form). Our main purpose here is to show
that the twistor spaces of self-dual Einstein $4$-manifolds of
negative scalar curvature admit such almost-K\"ahler structures.

   Recall that the twistor space of an oriented Riemannian $4$-manifold $M$ is the $2$-sphere bundle ${\cal Z}$ on $M$
consisting of the unit $(-1)$-%
eigenvectors of the Hodge star operator acting on $\Lambda^{2}TM$. The  $6$-manifold ${\cal Z}$ admits a natural
1-parameter family of Riemannian metrics $h_{t}$ such that the natural projection $\pi :{\cal Z} \to M$ is a
Riemannian submersion with totally geodesic fibres. These metrics are compatible with the almost-complex
structures $J_{1}$ and $J_{2}$ on ${\cal Z}$ introduced, respectively, by Atiyah, Hitchin \& Singer \cite{AHS} and
Eells \& Salamon \cite{ES}.

   In Section 3 we show that the first Chern form of the almost Hermitian manifold $({\cal Z},h_{t},J_2)$  always vanishes
which generalizes a result of Eells \& Salamon \cite{ES} stating that the almost-complex structure $J_{2}$ has
vanishing first Chern class. We obtain also an explicit formula for the first Chern form of $({\cal Z},h_{t},J_1)$
in terms of the curvature of the base manifold $M$. In the case when $M$ is self-dual the latter formula has been
actually given by P.Gauduchon \cite{G}.

   In Section 4 we obtain the precise geometric conditions on $M$ ensuring that the curvature of the Chern connection $D^n$
of $({\cal Z},h_{t},J_{n})$, $n=1,2$ is of type $(1,1)$.  Note that, in many cases,  this property of the
curvature simplifies the computation of the Chern numbers (cf. e.g. \cite{GBN}). We also study the problem when
the connection $D^n$  $n=1,2,$ has a constant holomorphic sectional curvature. The motivation behind is the open
question whether there are examples of non-K\"ahler Hermitian manifolds whose Chern connection is of non-zero
constant holomorphic sectional curvature (cf. \cite{Ba, BG}). Proposition~\ref{Proposition 5.5} shows that there
are no twistorial examples of such manifolds.

   In the last section we prove that the Nijenhuis tensor of the almost-complex structure $J_2$  is $D^2$-parallel
provided that the base manifold $M$ is Einstein and self-dual. Since in (real) dimension six, the Nijenhuis tensor
can be identified via the metric with a section of the canonical bundle, the result strengthens the fact that
$c_1({\cal Z}, J_2)=0$. If, in addition, $M$ is of negative scalar curvature $s$, then the twistor space $({\cal
Z},h_{t},J_{2})$, $t=-\displaystyle{\frac{12}{s}}$ is an almost- K\"ahler manifold with vanishing first Chern
form, the curvature of its Chern connection is of type $(1,1)$ and the Nijenhuis tensor of $J_2$ is parallel with
respect to it. Finally we note that the analogous statements for the twistor spaces of quaternionic-K\"ahler
manifolds are also valid .

\section{Preliminaries}

Let $M$ be a (connected) oriented Riemannian  $4$-manifold  with metric $g$. Then $g$ induces a metric on the
bundle $\Lambda ^{2}TM$ of $2$-vectors by the formula
$$
 g(X_{1}\land X_{2},X_{3}\land X_{4})=\frac 12[g(X_1,X_3)g(X_2,X_4)-g(X_1,X_4)g(X_2,X_3)]
$$
The Riemannian connection of $M$ determines  a  connection  on  the vector bundle  $\Lambda^{2}TM$ (both denoted
by $\nabla$) and the respective curvatures are related by
$$
 R(X,Y)(Z\land T)=R(X,Y)Z\land T + X\land R(Y,Z)T
$$
for $X,Y,Z,T\in \chi(M)$; $\chi(M)$ stands  for the Lie algebra  of  smooth vector fields on $M$. (For the
curvature tensor $R$  we  adopt the following definition $R(X,Y)=\nabla _{[X,Y]} - [\nabla _{X},\nabla _{Y}])$.
The  curvature operator ${\cal R}$ is the self-adjoint endomorphism of $\Lambda ^{2}TM$ defined by
$$
 g({\cal R}(X\land Y),Z\land T)=g(R(X,Y)Z,T)
$$
for all $X,Y,Z,T\in \chi (M)$. The Hodge star operator defines an endomorphism $*$ of $\Lambda^{2}TM$ with
$*^{2}=Id$. Hence
$$
 \Lambda^{2}TM=\Lambda^{2}_{+}TM\oplus \Lambda^{2}_{-}TM
$$
where $\Lambda^{2}_{\pm }TM$ are the subbundles of $\Lambda^{2}TM$ corresponding  to   the $(\pm 1)$-eigenvectors
of $*$. Let $(E_{1},E_{2},E_{3},E_{4})$ be  a  local  oriented orthonormal frame of $TM$. Set

\begin{equation}\label{eq 2.1}
\begin{array}{lll}
s_1=E_1\land E_2-E_3\land E_4&\ \ \ \ \ \ \ \ \ &
          \bar{s}_1=E_1\land E_2+E_3\land E_4 \\
s_2=E_1\land E_3-E_4\land E_2&\ \ \ \ \ \ \ \ \ &
          \bar{s}_2=E_1\land E_3+E_4\land E_2 \\
s_3=E_1\land E_4-E_2\land E_3&\ \ \ \ \ \ \ \ \ &
          \bar{s}_3=E_1\land E_4+E_2\land E_3
\end{array}
\end{equation}
Then $(s_{1},s_{2},s_{3})$ (resp.$(\bar{s}_{1},\bar{s}_{2},\bar{s}_{3})$ is  a  local  oriented  orthonormal frame
of $\Lambda ^{2}_{-}TM$ (resp.$\Lambda ^{2}_{+}TM$). The matrix of ${\cal R}$  with  respect  to  the frame
$(\bar{s}_{i},s_{i})$ of $\Lambda^{2}TM$ has the form
$$
{\cal R} =\left[
\begin{array}{cc}
A & B \\ {}^tB & C
\end{array} \right]
$$
where the $3\times 3$-matrices $A$ and $C$  are  symmetric  and  have  equal traces. Let ${\cal B}$, ${\cal
W}_{+}$ and ${\cal W}_{-}$ be the endomorphisms of $\Lambda ^{2}TM$ with matrices:
$$
{\cal B}=\left[
\begin{array}{cc}
 0 & B \\
{}^tB & 0
\end{array} \right], \ \
{\cal W}_{+}=\left[
\begin{array}{cc}
A-\displaystyle{\frac{s}{6}} I & 0 \\
0  & 0
\end{array}\right], \ \
{\cal W}_{-}=\left[
\begin{array}{cc}
0 & 0 \\
0&C-\displaystyle{\frac{s}{6}} I
\end{array} \right]
$$
where $s$ is the scalar curvature and $I$ is the unit $3\times 3$-matrix. Then
$$
{\cal R}=\frac{s}{6} Id+{\cal B}+{\cal W}_{+}+{\cal W}_{-}
$$
is the irreducible decomposition of ${\cal R}$ under the action of $SO(4)$  found by Singer \& Thorpe \cite{ST}.
Note that ${\cal B}$ and ${\cal W}={\cal W}_{+}+{\cal W}_{-}$ represent the traceless Ricci  tensor  and  the
Weyl  conformal  tensor, respectively. The manifold $M$ is called self-dual (anti-self-dual) if ${\cal W}_{-}=0$
(${\cal W}_{+}=0$). It is Einstein exactly when ${\cal B}=0$.

The twistor space of $M$ is the subbundle ${\cal Z}$ of $\Lambda^2_-TM$ consisting of all unit vectors. The
Riemannian connection $\nabla$ of $M$ gives rise to a splitting $T{\cal Z}={\cal H} \oplus {\cal V}$ of the
tangent bundle of ${\cal Z}$ into horizontal and vertical components. More precisely, let $\pi :\Lambda^2_-TM\to
M$ be the natural projection. By definition, the vertical space at $\sigma\in{\cal Z}$ is ${\cal V}_{\sigma
}=$Ker$\pi_{*\sigma}$ ($T_{\sigma}{\cal Z}$ is always considered as a subspace of $T_{\sigma}(\Lambda^2_-TM)$).
Note that ${\cal V}_{\sigma }$ consists of those vectors of $T_{\sigma }{\cal Z}$ which are tangent to the fibre
${\cal Z}_{p}=\pi^{-1}(p)\cap {\cal Z}$, $p=\pi(\sigma)$, of ${\cal Z}$ through the point $\sigma$. Since ${\cal
Z}_{p}$ is the unit sphere in the vector space $\Lambda^2_{-}T_pM$, ${\cal V}_{\sigma}$ is the orthogonal
complement of $\sigma$ in $\Lambda^2_{-}T_{p}M$. Let $\xi$ be a local section of ${\cal Z}$ such that
$\xi(p)=\sigma $. Since $\xi$ has a constant length, $\nabla_{X}\xi\in {\cal V}_{\sigma}$ for all $X\in T_{p}M$.
Given $X\in T_{p}M$, the vector $X^{h}_{\sigma}=\xi_{*}X-\nabla_{X}\xi\in T_\sigma{\cal Z}$ depends only on $p$
and $\sigma $. By definition, the horizontal space at $\sigma $ is ${\cal H}_{\sigma}=\{X^{h}_{\sigma}:X\in
T_{p}M\}$. Note that the map $X \to X^{h}_{\sigma }$ is an isomorphism between $T_{p}M$ and ${\cal H}_{\sigma }$
with inverse map $\pi_*\mid {\cal H}_{\sigma }$.

Let $(U,x_1,x_2,x_3,x_4)$ be a local coordinate system of  $M$  and let $(E_1,E_2,E_3,E_4)$ be an oriented
orthonormal frame of $TM$ on $U$. If $(s_1,s_2,s_3)$ is the local frame of $\Lambda^{2}_{-}TM$ defined by
(\ref{eq 2.1}),  then $\tilde x_i=x_i \circ\pi$, $y_j(\sigma)=g(\sigma,(s_j\circ\pi)(\sigma))$, $1\le i\le 4$,
$1\le j\le 3$, are local  coordinates of $\Lambda^2_-TM$ on $\pi^{-1}(U)$. For each vector field
$$
 X=\sum^{4}_{i=1}X^{i}\frac{\partial }{\partial x_i}
$$
on $U$ the horizontal lift $X^{h}$ of $X$ on $\pi^{-1}(U)$ is given by
\begin{equation}\label{eq 2.2}
X^{h}=\sum^{4}_{i=1}(X^i\circ\pi)\frac{\partial}{\partial\tilde x_i} -
\sum^{3}_{j,k=1}y_{j}g(\nabla_Xs_j,s_k)\circ\pi\frac{\partial}{\partial y_k}.
\end{equation}

Let $\sigma\in{\cal Z}$ and $\pi(\sigma)=p$. Using (\ref{eq 2.2}) and  the standard  identification
$T_{\sigma}(\Lambda^2_-T_pM) \\
\cong \Lambda^2_-T_pM$ one gets that
\begin{equation}\label{eq 2.3}
[X^h,Y^h]_{\sigma}-[X,Y]^h_{\sigma}=R_p(X\land Y)\sigma
\end{equation}
for all $X,Y\in \chi(U)$.

Each  point  $\sigma\in{\cal Z}$ defines a complex structure $K_{\sigma}$ on $T_pM$ by
\begin{equation}\label{eq 2.4}
g(K_{\sigma}X,Y)=2g(\sigma,X\land Y), \ \ X,Y\in T_{p}M.
\end{equation}
Note that $K_{\sigma }$ is compatible with the metric  $g$  and  the opposite orientation of $M$ at $p$. The
2-vector $2\sigma$ is dual  to  the fundamental 2-form of $K_{\sigma }$.

Denote  by  $\times$  the  usual  vector  product  in  the  oriented 3-dimensional vector space
$\Lambda^{2}_{-}T_{p}M$, $p\in M$. Then it is easily checked that
\begin{equation}\label{eq 2.5}
g(R(a)b,c)=-g({\cal R}(a),b\times c))
\end{equation}
for $a\in\Lambda^2T_{p}M$, $b,c\in \Lambda ^{2}_{-}T_{p}M$  and
\begin{equation}\label{eq 2.6}
g(\sigma\times V,X\land K_{\sigma }Y)=%
g(\sigma\times V,K_{\sigma }X\land Y)=-g(V,X\land Y)
\end{equation}
for  $V\in {\cal V}_{\sigma}$, $X,Y\in T_{p}M$.

   It is also easy to check that for any $\sigma, \tau\in{\cal Z}$ with $\pi(\sigma)=\pi(\tau)$ we have
\begin{equation}\label{eq 2.6C}
K_{\sigma}\circ K_{\tau}=-g(\sigma,\tau)Id-K_{\sigma\times\tau}
\end{equation}

  Following \cite{AHS} and \cite{ES} define two almost--complex  structures $J_{1}$ and $J_{2}$ on ${\cal Z}$ by
$$
J_{n}V=(-1)^{n} \sigma\times V \mbox{ for } \ V\in {\cal V}_{\sigma }
$$
$$
J_{n}X^{h}_{\sigma }=(K_{\sigma}X)^{h}_{\sigma} \mbox{ for } \ X\in T_{p}M, p=\pi (\sigma ).
$$
It is well-known \cite{AHS} that $J_{1}$ is integrable (i.e. comes  from  a complex  structure) if and only if $M$
is  self-dual.   Unlike   $J_{1}$,   the almost-complex structure $J_{2}$ is never integrable \cite{ES}.

   Let $h_{t}$ be the Riemannian metric on ${\cal Z}$ given by
$$
h_{t}=\pi ^{*}g+tg^{v}
$$
where $t>0$, $g$ is the metric of $M$ and $g^{v}$ is  the  restriction  of the metric of $\Lambda^{2}TM$ on the
vertical distribution ${\cal V}$. Then $\pi:(Z,h_t)\to (M,g)$ is a Riemannian submersion with totally geodesic
fibres and the almost-complex structures $J_{1}$ and $J_{2}$ are compatible with the metrics $h_{t}$. Denote by $D
(=D_t)$ the Levi-Civita connection of $({\cal Z},h_t)$. Let $\sigma$ be a point of ${\cal Z}$, $X,Y$ vector fields
on $M$ near the point $\pi(\sigma)$ and $A$ a vertical vector field near $\sigma$. It is not hard to see
(cf.\,e.g.\,\cite{DM2}) that
\begin{equation}\label{eq 2.7}
(D_{X^h}Y^h)_\sigma=(\nabla_XY)^h_\sigma+\frac 12 R(X\land Y)\sigma
\end{equation}
\begin{equation}\label{eq 2.8}
(D_AX^h)_\sigma={\cal H}(D_{X^h}A)_\sigma=\frac t2 (R(\sigma\times A)X)^h_\sigma
\end{equation}

\section{The first Chern forms of twistor spaces}

 Given an almost-Hermitian manifold $(N,g,J)$, denote by $\nabla $ the Levi-Civita connection of $g$. Then the Chern
connection $\tilde{\nabla }$  of $(N,g,J)$ is defined by (cf. e.g. \cite[Th.6.1]{GBN}):
\begin{equation}\label{eq 5.1}
\begin{array}{c}
g(\tilde{\nabla}_XY,Z)=g(\nabla_XY,Z)+\frac 12g((\nabla_XJ)(JY),Z)\\
+\frac 14g((\nabla_ZJ)(JY)-(\nabla_YJ)(JZ)-(\nabla_{JZ}J)(Y)+(\nabla_{JY}J)(Z),X)
\end{array}
\end{equation}

It is one of the distinguished 1-parameter family of Hermitian
connections defined by P.Gauduchon \cite{Gaud}:

\begin{equation}
\begin{array}{c}
g(\tilde{\nabla}^u_XY,Z)=g(\nabla_XY,Z)+\frac 12g((\nabla_XJ)(JY),Z)\\
+\frac
{u}{4}g((\nabla_ZJ)(JY)-(\nabla_YJ)(JZ)-(\nabla_{JZ}J)(Y)+(\nabla_{JY}J)(Z),X)
\end{array}
\end{equation}

The Chern connection corresponds to $u=1$. Let
$\Omega(X,Y)=g(JX,Y)$ be the K\"ahler form of $(N,g,J)$ and
$\delta\Omega$ the codifferential of $\Omega$ with respect to
$\nabla$. Denote by $\varphi$ and $\psi$ the 2-forms on $N$
defined by
\begin{equation}\label{eq 5.2}
\varphi(X,Y)=Trace(Z\to g((\nabla_XJ)(JZ),(\nabla _YJ)(Z)))
\end{equation}
\begin{equation}\label{eq 5.3}
\psi(X,Y)=\rho^*(X,JY)
\end{equation}
where $\rho^*$ is the $*$-Ricci tensor of $(N,g,J)$. Recall that $\rho^*$ is given by
$$\rho^*(X,Y)=Trace(Z\to R(JZ,X)JY),$$
where $R$ is the curvature tensor of $\nabla$. The formula in the
next Lemma appears in \cite{Gaud} without proof, so for sake of
completeness we provide its proof here.

\begin{lem}\label{Lemma 5.1} The first Chern form $\gamma^u$ of the connection $\tilde{\nabla}^u$ on an almost Hermitian manifold
$(N,g,J)$ is given by
$$
 8\pi\gamma^u=-\varphi-4\psi+2ud\delta\Omega
$$
\end{lem}

\noindent {\it Proof.}
Denote by $\hat\nabla$ the connection on $N$ defined by
$$
\hat\nabla_{X}Y=\nabla_XY+\frac 12(\nabla _XJ)(JY), \ \ X,Y\in\chi(N).
$$
Note that $\hat\nabla g=0$ and $\hat\nabla J=0$. Let $S$ be the
(1,2)-tensor field on $N$
defined by
\begin{equation}\label{eq 5.4}
g(S(X,Y),Z)=\frac 14g((\nabla_ZJ)(JY)-(\nabla_YJ)(JZ)-%
(\nabla_{JZ}J)(Y)+(\nabla_{JY}J)(Z),X)
\end{equation}
Then
$$
\tilde{\nabla}^u_XY=\hat\nabla_XY+uS(X,Y).
$$
Bellow we consider only the case $u=1$ since the general case
follows immediately from it. It is easy to check that $S$ has the
following properties:
\begin{equation}\label{eq 5.5}
g(S(X,Y),Z)=-g(S(X,Z),Y)
\end{equation}
\begin{equation}\label{eq 5.6}
S(X,JY)=JS(X,Y)
\end{equation}
\begin{equation}\label{eq 5.7}
S(X,JX)=-S(JX,X)=\frac 14((\nabla_XJ)(X)+(\nabla_{JX}J)(JX))
\end{equation}
\begin{equation}\label{eq 5.8}
g((\hat\nabla_YS)(JX,JX),X)=0.
\end{equation}
 A straightforward computation shows that the curvature tensors $R$, $\hat{R}$ and $\tilde{R}$ of $\nabla$,
$\hat\nabla$  and $\tilde{\nabla}$ are related by
\begin{equation}\label{eq 5.9}
\begin{array}{c}
4\hat{R}(X,Y,Z,W)=2R(X,Y,Z,W)+2R(X,Y,JZ,JW)\\
+g((\nabla_XJ)(Z),(\nabla_YJ)(W))-g((\nabla _XJ)(W),(\nabla_YJ)(Z))
\end{array}
\end{equation}
\begin{equation}\label{eq 5.10}
\begin{array}{c}
\tilde{R}(X,Y,Z,W)=\hat{R}(X,Y,Z,W)-g((\hat\nabla_XS)(Y,Z),W)+%
g((\hat\nabla_YS)(X,Z),W)\\
+g(S(X,W),S(Y,Z))-g(S(Y,W),S(X,Z))-g(S(\hat{T}(X,Y),Z),W)
\end{array}
\end{equation}
where $\hat{T}$ is the torsion of $\hat\nabla$.

  Now fix a point $p\in N$ and choose an orthonormal frame $E_1,\ldots,E_n,JE_1,\ldots,JE_n$ near $p$ such that
$\hat\nabla E_i\mid_{p}=0$, $i=1,\ldots,n$. Then, using (\ref{eq 5.5}), (\ref{eq 5.6}) and (\ref{eq 5.10}), one
gets:
$$
4\pi\gamma(X,Y)=\sum^{2n}_{k=1}\tilde R(X,Y,E_k,JE_k)=%
\sum^{2n}_{k=1}[\hat{R}(X,Y,E_k,JE_k)
$$
$$
- X(g(S(Y,E_k),JE_k))+Y(g(S(X,E_k),JE_k))+g(S([X,Y,],E_k),JE_k)]
$$
at the point $p$.

Formula (\ref{eq 5.9}) together with the first Bianchi identity gives:
$$
\sum^{2n}_{k=1}\hat{R}(X,Y,E_k,JE_k)=-2\psi(X,Y)-\frac 12\varphi(X,Y)
$$

Moreover, by (\ref{eq 5.4}), one has:
$$
\sum^{2n}_{k=1}g(S(X,E_k),JE_k)=-\delta\Omega(X)
$$
and  the lemma follows from the above identities.

 Now let $M$ be an oriented Riemannian 4-manifold with twistor
space ${\cal Z}$. Let $D (=D_{t})$ be the Levi-Civita connection of $({\cal Z},h_{t})$. Denote by $D^{n}
(=D^{n}_{t})$ the Chern connection of the almost-Hermitian manifold $({\cal Z}$, $h_{t}$, $J_{n})$, $n=1,2$ and by
$\gamma_{t,n}$ its first Chern form. In the case when the base manifold $M$ is self-dual an explicit formula for
the first Chern form of $D^1$ has been given by P.Gauduchon \cite{G}. Here we compute the first Chern forms
$\gamma_{t,n}$, $n=1,2$, of the twistor space of an arbitrary oriented 4-manifold $M$. To do this we shall use the
following formulas for the covariant derivative of the almost-complex structure $J_{n}$ with respect to the
Levi-Civita connection $D$ (\cite{M}):

\begin{lem}\label{Lemma 5.2} For any $\sigma\in{\cal Z}$, $A\in{\cal V}_{\sigma}$  and $X,Y\in T_pM$,
$p=\pi(\sigma)$, one has
$$
h_t((D_{X^h}J_n)(Y^h),A)=\frac t2 [(-1)^ng({\cal R}(A),X\land Y)-%
g({\cal R}(\sigma\times A),X\land K_{\sigma}Y)]
$$
$$
h_t(D_AJ_n)(X^h),Y^h)=\frac t2 g({\cal R}(\sigma\times A),X\land K_{\sigma}Y+%
K_{\sigma}X\land Y)+2g(A,X\land Y).
$$
where $K_{\sigma}$ is the complex structure on $T_pM$ defined via {\rm (\ref{eq 2.4})}. Moreover,
$$
h_t(D_EJ_n)(F),G)=0
$$
whenever $E,F,G$ are horizontal vectors or at least two of them are vertical vectors.
\end{lem}

 We shall also need the following formula for the $*$-Ricci tensor $\rho^*_{t,n}$ of $({\cal Z}$, $h_{t}$, $J_{n})$
\cite{DMG2}:

\begin{lem}\label{Lemma 5.3} Let $E,F\in T_{\sigma }{\cal Z}$
 and $X=\pi _{*}E$ ,$Y=\pi _{*}F$ ,$A={\cal V}E$, $B={\cal V}F$. Then
$$
\begin{array}{ll}
\rho^*_{t,n}(E,F)&= [1 +(-1)^{n+1}]g({\cal R}(\sigma ),%
X\land K_{\sigma }Y)-\frac t2 g(R(X\land K_{\sigma }Y)\sigma,R(\sigma)\sigma) \\
&+ \frac t4\mbox{ Trace}(Z\to g(R(X\land Z)\sigma,R(K_{\sigma }Z\land
K_{\sigma}Y)\sigma))\\
&+ \frac t4 (-1)^{n+1}\mbox{ Trace}({\cal V}_{\sigma}\ni C\to g(R(C)X,R(\sigma\times
C)K_{\sigma }Y))\\
&+ \frac t2 (-1)^{n}g((\nabla _{X}{\cal R})(\sigma),B)+\frac t2 g((\nabla _{K_{\sigma }Y}{\cal
R})(\sigma ),\sigma\times A)\\ &+ [1 +(-1)^{n+1}tg({\cal R}(\sigma ),\sigma )]g(A,B)\\
&+ (-1)^{n+1}\frac{t^2}4\mbox{ Trace}(Z\to g(R(\sigma\times A)K_{\sigma }Z,R(B)Z)).
\end{array}
$$
%where $K_{\sigma }$ is the complex structure on $T_{p}M$, $p=\pi(\sigma)$, determined  by $\sigma$ via {\rm
%(\ref{eq 2.4}}.
\end{lem}

  Now we are ready to prove the following

\begin{prop}\label{Proposition 5.3} The first Chern form $\gamma_{t,n}$ of the twistor space $({\cal
Z},h_t,J_n)$, $n=1,2,$ is given by
$$
2\pi\gamma_{t,n}(E,F)=[1+(-1)^{n+1}][g({\cal R}(\sigma),X\land Y)+%
g(A,\sigma\times B)]
$$
where $E,F\in T_{\sigma}{\cal Z}$  and $X=\pi_{*}E$, $Y=\pi_{*}F$, $A={\cal V}E$, $B={\cal V}F$.
\end{prop}

\noindent {\it Proof.}  Denote by $\varphi_{t,n}$ and $\psi_{t,n}$ the 2-forms on ${\cal Z}$ defined by (\ref{eq
5.2}) and (\ref{eq 5.3}), respectively. Let $\Omega_{t,n}$ be the K\"ahler form of $({\cal Z},h_t,J_n)$, $n=1,2$.
By Lemma~\ref{Lemma 5.1}, we have
$$8\pi\gamma_{t,n}=-\varphi_{t,n}-4\psi_{t,n}+2d\delta\Omega_{t,n}.$$
Let $U$ be an $h_t$-unit vertical vector at $\sigma$. Then, using Lemma~\ref{Lemma 5.2}, (\ref{eq 2.5}) and
(\ref{eq 2.6}), one gets:
$$
\begin{array}{rl}
\varphi_{t,n}(E,F)&
=tTrace(Z\to g(R(X\land Z)\sigma,R(K_{\sigma }Z\land Y)\sigma ))\\
&
+(-1)^{n+1}tTrace({\cal V}_{\sigma}\ni C\to g(R(C)X,R(\sigma\times C)Y))\\
&
+t^2Trace (Z\to g(R(\sigma\times A)K_{\sigma }Z,R(\sigma\times B)Z))\\
&
-2tg({\cal R}(\sigma\times A),B)+2tg({\cal R}(\sigma\times B),A)+%
4g(\sigma\times A,B)
\end{array}
$$
Since $\psi_{t,n}(E,F)=c^*_{t,n}(E,J_nF)$, it follows from Lemma~\ref{Lemma 5.3} that
\begin{equation}\label{eq 5.11}
\begin{array}{c}
4\pi\gamma_{t,n}(E,F)=2[1+(-1)^{n+1}][g({\cal R}(\sigma),X\land Y)+%
g(A,\sigma\times B)]\\
+t[g({\cal R}(\sigma\times A),B)-g({\cal R}(\sigma\times B),A)+%
2g({\cal R}(\sigma),\sigma)g(A,\sigma\times B)\\
-g((\nabla_X{\cal R})(\sigma),\sigma\times B)+%
g((\nabla_Y{\cal R})(\sigma),\sigma\times A)-g(R(X\land Y)\sigma,%
R(\sigma)\sigma)] \\
+d\delta\Omega _{t,n}(E,F).
\end{array}
\end{equation}
It is easy to check by means of Lemma~\ref{Lemma 5.2} and the identity (\ref{eq 2.5}) that the 1-form $\omega=-
1/t\delta\Omega_{t,n}$ is given by $\omega(E)=g({\cal V}E,R(\sigma)\sigma )$ for $E\in T_{\sigma}{\cal Z}$. Next
we shall compute the differential of the form $\omega$. Since $\sigma \to R(\sigma)\sigma$ is a vertical vector
field on ${\cal Z}$, one has by (\ref{eq 2.3}):
\begin{equation}\label{eq 5.12}
(d\omega)_{\sigma}(X^h,Y^h)=-\omega_{\sigma}([X^h,Y^h])=%
-g(R(X\land Y)\sigma,R(\sigma)\sigma); \ \  X,Y\in \chi(M)
\end{equation}
 Now let $s$ be a local section of ${\cal Z}$ such that $s(p)=\sigma$ and $\nabla s\mid_p=0.$ If $B$ is a
vertical vector field on ${\cal Z}$ and $X$ is a vector field on $M$, it follows easily from (\ref{eq 2.2}) that
$$
[X^h,B]_{\sigma}=\nabla_{X_p}(B\circ s)
$$
where $B\circ s$ is considered as a section of $\Lambda^2_-TM$. Then
$$
(d\omega)_\sigma(X^h,B)=s_*(X_p)(\omega(B))-\omega_{\sigma}([X^h,B])
$$
$$
=X_p(g(B\circ s,R(s)s))-g(\nabla_{X_p}(B\circ s),R(\sigma)\sigma)
$$
and, using (\ref{eq 2.5}), one gets:
\begin{equation}\label{eq 5.13}
(d\omega)_{\sigma}(X^h,B)=-g((\nabla_{X_p}{\cal R})(\sigma),%
\sigma\times B_{\sigma})
\end{equation}
 Finally, we will show that
\begin{equation}\label{eq 5.14}
(d\omega)_{\sigma}(A,B)=g({\cal R}(\sigma\times A),B)-%
g({\cal R}(\sigma\times B),A)+2g({\cal R}(\sigma),\sigma)g(A,\sigma\times B)
\end{equation}
for any vertical vectors $A$ and $B$ at $\sigma $.

Let $(s_1,s_2,s_3)$ be a local frame of $\Lambda^{2}_{-}TM$ defined by (\ref{eq 2.1}) such that $s_{1}(p)=\sigma$
and let $y_j(\tau)=g(\tau,(s_j\circ\pi)(\tau))$, $\tau\in\Lambda^2_-TM$, $1\le j \le 3$. Set ${\displaystyle
U=-y_2\frac{\partial}{\partial y_1}+y_1\frac{\partial}{\partial y_2} }$. Then ${\displaystyle
J_1U=y_1y_3\frac{\partial}{\partial y_1}+%
y_2y_3\frac{\partial}{\partial y_2}}$ ${\displaystyle -(1-y^2_3)\frac{\partial}{\partial y_3}}$ and $(U,J_1U)$ is
a local frame of the vertical bundle ${\cal V}$ near the point $\sigma$ such that $[U,J_1U]_{\sigma }=0$. It is
enough to check (\ref{eq 5.14}) for $A=U_{\sigma }$ and $B=J_1U_{\sigma }$. Using (\ref{eq 2.5}), one gets:
$$
\omega(U)=\sum^{3}_{j=1}y_j[y_1y_3g({\cal R}(s_j),s_1)\circ\pi+%
y_2y_3g({\cal R}(s_j),s_2)\circ\pi- (1-y^2_3)g({\cal R}(s_{j}),s_{3})\circ\pi]
$$
$$
\omega(J_1U)=\sum^{3}_{j=1}y_j[y_2g({\cal R}(s_j),s_1)\circ\pi-%
y_1g({\cal R}(s_j),s_2)\circ\pi]
$$
Then
$$
(d\omega)_\sigma(U,J_1U)=(\frac{\partial}{\partial y_2})_\sigma(\omega(J_1U))+%
(\frac{\partial }{\partial y_3})_{\sigma}(\omega(U))
$$
$$
=-g({\cal R}(s_3),s_3)-g({\cal R}(s_2),s_{2})+2g({\cal R}(s_1),s_1)
$$
which proves (\ref{eq 5.14}).

Now the proposition follows from (\ref{eq 5.11}) - (\ref{eq 5.14}).

\section{Curvature properties of the Chern connection of twistor spaces}

  In this section we consider the problem when the curvature tensor $R_{t,n}$
of the Chern connection $D^n$ of $({\cal Z},h_t,J_n)$, $n=1,2,$ is of type $(1,1)$, i.e.
$R_{t,n}(J_nE,J_nF)G=R_{t,n}(E,F)G$ for all $E,F,G\in T{\cal Z}$. We also study the problem when this connection
has a constant holomorphic sectional curvature.

\begin{prop}\label{Proposition 5.4}
\begin{enumerate}
\item[$(i)$] The curvature tensor $R_{t,1}$ is of type $(1,1)$ if and only if the base manifold $M$
is self-dual.
\item[$(ii)$] The curvature tensor $R_{t,2}$ is of type $(1,1)$  if and only if the base manifold $M$ is Einstein and
self-dual.
\end{enumerate}
\end{prop}

\noindent {\it Proof.} (i) If $R_{t,1}$ is of type $(1,1)$, then $\gamma _{t,1}$ is an (1,1)-form with respect to
$J_{1}$. This together with Propostion~\ref{Proposition 5.3} gives:
$$
g({\cal R}(\sigma),X\land Y-K_{\sigma}X\land K_{\sigma}Y)=0
$$
for all $\sigma\in{\cal Z}$ and $X,Y\in T_pM$, $p=\pi(\sigma )$. Since the 2-vectors of the form $X\land Y -
K_\sigma X\land K_{\sigma }Y$ span the vertical space at $\sigma $, it is easy to see that the latter identity
implies the self-duality of $M$.

Conversely, if $M$ is self-dual, the almost-complex structure $J_{1}$ is integrable \cite{AHS} and, as it is
well-known (cf. e.g. \cite[Lemma 2.1]{GBN}), the curvature of the Chern connection $D^{1}$ is of type $(1,1)$.

(ii) Given a point $\sigma \in {\cal Z}$ and $X,Y\in T_{p}M$, $p=\pi (\sigma )$, denote by $A(X,Y)$ and $B(X,Y)$
the vertical vectors at $\sigma $ defined by
$$
A(X,Y)=\frac 14R(X\land Y+K_{\sigma }X\land K_{\sigma }Y)\sigma-%
\frac 1{4t} (K_{\sigma }X\land Y+X\land K_{\sigma }Y)
$$
$$
B(X,Y)=\frac 14\sigma\times R(K_{\sigma }X\land Y-%
X\land K_{\sigma }Y)\sigma+\frac 1{4t}(K_{\sigma }X\land Y+%
X\land K_{\sigma }Y)
$$
where $K_{\sigma }$ is the complex structure on $T_{p}M$ corresponding to $\sigma $ via (\ref{eq 2.4}). Using
Lemma~\ref{Lemma 5.2}, formulas (\ref{eq 2.5}),(\ref{eq 2.6}), (\ref{eq 2.7}) and the identity
$\sigma\times (K_{\sigma }X\land Y+X\land K_{\sigma }Y)=%
X\land Y-K_{\sigma }X\land K_{\sigma }Y$ one gets from (\ref{eq 5.1}) that
$$
D^{2}_{X^{h}}Y^{h} = (\nabla _{X}Y)^{h} + A(X,Y) + B(X,Y)
$$
$$
D^{1}_{X^{h}}Y^{h} = D^{2}_{X^{h}}Y^{h} - 2B(X,Y)
$$
$$
D^{1}_{V}X^{h} = D^{2}_{V}X^{h}
$$
for $X,Y\in \chi (M)$ and $V\in {\cal V}$.

Now one obtains easily that
%\begin{equation}\label{eq 5.15}
$$
\begin{array}{rl}
R_{t,1}(X^{h},Y^{h},Z^{h},T^{h})%
&=R_{t,2}(X^{h},Y^{h},Z^{h},T^{h})\\
&+2t[g(A(X,Z),B(Y,T))+g(A(Y,T),B(X,Z))\\
&-g(A(Y,Z),B(X,T) -g(A(X,T),B(Y,Z))]\\
\end{array}
%\end{equation}
$$
and

$$
g(A(X,Z),B(Y,T)) - g(A(K_{\sigma }Y,T),B(K_{\sigma }X,Z))
$$
$$
=\frac 1{8t}[g({\cal R}(X\land K_{\sigma }Z+K_{\sigma}X\land Z),%
K_{\sigma }Y\land T-Y\land K_{\sigma }T)
$$
$$
-g({\cal R}(X\land Z+K_{\sigma }X\land K_{\sigma }Z),Y\land T-%
K_{\sigma }Y\land K_{\sigma }T)]
$$
Theses formulas together with the first Bianchi identity give:
\begin{equation}\label{eq 5.16}
\begin{array}{l}
\ \ R_{t,1}(X^h,Y^h,Z^h,T^h)-R_{t,1}(J_1X^h,J_1Y^h,Z^h,T^h)\\
=R_{t,2}(X^h,Y^h,Z^h,T^h)-R_{t,2}(J_2X^h,J_2Y^h,Z^h,T^h)\\
-\frac 12 g({\cal R}(X\land Y-K_{\sigma }X\land K_{\sigma }Y,Z\land T+%
K_{\sigma }Z\land K_{\sigma }T)
\end{array}
\end{equation}
 A similar computation gives
\begin{equation}\label{eq 5.17}
\begin{array}{l}
\ \ R_{t,1}(X^h,Y^h,V,J_1V)-R_{t,1}(J_1X^h,J_1Y^h,V,J_1V)\\
=-R_{t,2}(X^h,Y^h,V,J_2V)+R_{t,2}(J_2X^h,J_2Y^h,V,J_2V) \\ %
\hspace{0.4cm}+ g({\cal R}(X\land Y-K_{\sigma }X\land K_{\sigma}Y),\sigma )
\end{array}
\end{equation}
for any $h_{t}$-unit vertical vector $V$ at $\sigma $.

Now assume that $R_{t,2}$ is of type $(1,1)$ with respect to $J_{2}$. Then it follows from (\ref{eq 5.16}) and
(\ref{eq 5.17}) that
$$
\gamma _{t,1}(X^h,Y^h)-\gamma _{t,1}(J_1X^h,J_1Y^h)=0
$$
which implies, as we have seen in the proof of (i), that $M$ is self-dual. Hence, by (i), $R_{t,1}$ is of type
(1,1) with respect to $J_{1}$ and the identity (\ref{eq 5.16}) becomes
$$
g({\cal R}(X\land Y-K_{\sigma }X\land K_{\sigma}Y,Z\land T+%
K_{\sigma }Z\land K_{\sigma }T)=0
$$
for $X,Y,Z,T\in\chi(M)$. Since $M$ is self-dual, this implies ${\cal B}=0$, i.e. $M$ is Einstein.

Conversely, let $M$ be Einstein and self-dual. Then the
almost-Hermitian manifold $({\cal Z},h_t,J_2)$ is quasi-K\"ahler
\cite{M}. On the other hand, according to \cite[Theorem,
(i)]{DMG}, its Riemannian curvature tensor satisfies the identity
$R(E,F,G,H)=\\R(JE,JF,JG,JH)$. Now it follows from
\cite[Th.6.2(ii)]{GBN} that the curvature tensor $R_{t,2}$  of the
Chern connection $D^{2}$ is of type $(1,1)$.

  Next, we study the problem when the Chern connections $D^{1}$ and $D^{2}$ of a twistor space have constant
holomorphic sectional curvatures.

\begin{prop}\label{Proposition 5.5} The Chern connection $D^{1}$ of the almost-Hermitian manifold \\
$({\cal Z},h_{t},J_{1})$ has a constant holomorphic sectional curvature $\kappa$ if and only if $\kappa>0$, the
base manifold  $M$ is of constant sectional curvature $\kappa$ and $t=1/\kappa$.

The holomorphic sectional curvature of the Chern connection $D^{2}$ of $({\cal Z},h_{t},J_{2})$ is never constant.
\end{prop}

\noindent {\it Proof.}  Let us note that if $(N,g,J)$ is an almost-Hermitian manifold, then the holomorphic
sectional curvatures $H$ and $\tilde H$ of the Levi-Civita connection $\nabla$ and the Chern connection
$\tilde\nabla$ are related by
\begin{equation}\label{eq H}
\begin{array}{l}
\tilde{H}(X)=H(X)+\displaystyle{\frac 18}(\|(\nabla_XJ)(X)\|^2+\|(\nabla_{JX}J)(JX)\|^2)\\
                                                                          \\
\hspace{2.8cm} +\displaystyle{\frac 34} g((\nabla_XJ)(X),(\nabla_{JX}J)(JX)).
\end{array}
\end{equation}
This easily follows from (\ref{eq 5.4}) - \ref{eq 5.9}).

  Denote by $\tilde{H}_{t,n}$ the holomorphic sectional curvature of the Chern connection $D^{n}$ of
$({\cal Z},h_{t},J_{n})$, $n=1,2$. Then using the explicit formula for the sectional curvature of $({\cal Z},h_t)$
given in \cite[Proposition 3.5]{DM2}, formula (\ref{eq H}) and Lemma~\ref{Lemma 5.2}, we obtain:
$$
\tilde{H}_{t,n}(X^h)=R(X,K_\sigma X,X,K_\sigma X)-%
\frac t2\|R(X\land K_\sigma X)\sigma\|^{2}_{g}
$$
where $K_{\sigma}$ is the complex structure on $T_{p}M$, $p=\pi(\sigma)$, defined by (\ref{eq 2.4}).

Assume that $\tilde{H}_{t,n}\equiv \kappa$. Then, for every $\sigma \in {\cal Z}$ and $X\in T_{p}M$, $p=\pi(\sigma
)$, $\|X\|=1$, one has:
\begin{equation}\label{eq 5.18}
\kappa =R(X,K_\sigma X,X,K_\sigma X)-\frac t2\|R(X\land K_\sigma X)\sigma\|^2_g
\end{equation}

 Let $s_1,s_2,s_3$ be local sections of ${\cal Z}$ defined by (\ref{eq 2.1}) and let
$${\displaystyle\sigma=\sum^{3}_{i=1}\lambda_is_i},\hspace{0.2cm} {\displaystyle\sum^{3}_{i=1}\lambda^2_i=1}.$$
Denote by $K_i$ the complex structure on $T_{p}M$ determined by $s_{i}(p)$ and set
$$
a_{ij}=g({\cal R}(s_i),X\land K_jX), \ \ b_{ij}=g({\cal R}(X\land K_iX),X\land K_jX)
$$
Then
$$
\|R(X\land K_\sigma X)\sigma\|^2_g=%
\sum^{3}_{i=1}g({\cal R}(\sigma\times s_i), \ \
X\land K_\sigma X)^2=\sum^{3}_{i=1}\Bigl(\sum^{3}_{j=1}\lambda_ja_{ij}\Bigr)^2-%
(\sum^{3}_{i,j}\lambda_i\lambda_ja_{ij} )^2
$$
and
$$
R(X,K_\sigma X,X,K_\sigma X)=\sum^{3}_{i,j=1}\lambda_{i}\lambda _{j}b_{ij}
$$
Varying $(\lambda_1,\lambda_2,\lambda_3)$ over the unit sphere $S^2$, one gets from (\ref{eq 5.18})that
$$
b_{ii}-\frac t2\sum^{3}_{k=1}a^2_{ki}+\frac t2a^2_{ii}=\kappa
$$
$$
b_{ii}+b_{jj}-\frac t2\sum^3_{k=1}(a^2_{ki}+a^2_{kj})+%
\frac t2(a_{ij}+a_{ji})^2+ta_{ii}a_{jj}=2\kappa
$$
$$
b_{ij}+b_{ji}-t\sum^{3}_{k=1}a_{ki}a_{kj}+%
ta_{ii}(a_{ij}+a_{ji})=0
$$
for $1\le i\neq j\le 3$. These identities imply $a_{ii}=a_{jj}$ and $a_{ij}=-a_{ji}$ for $i\neq j$, i.e.
$$
g({\cal R}(s_i),X\land K_iX) =g({\cal R}(s_j),X\land K_jX)
$$
$$
g({\cal R}(s_i),X\land K_jX)=-g({\cal R}(s_j),X\land K_iX), \ \ \ i\neq j
$$
Now varying $X$ over the unit sphere of $T_{p}M$ gives:
$$
g({\cal R}(s_{i}),s_{j})=\delta _{ij}g({\cal R}(s_{1}),s_{1})
$$
$$
g({\cal R}(s_{i}),\bar{s}_{j})=0, \ \ \  1\le i,j\le 3
$$
Hence $M$ is Einstein and self-dual. Since $X\land K_\sigma X\in \R.\sigma%
\oplus\Lambda_+^2T_pM$ for any $X\in T_pM$, it follows that $R(X\land K_\sigma X)\sigma=0$ and (\ref{eq 5.18})
shows that $M$ is of constant sectional curvature equal to $\kappa$. In this case one obtains easily from
\cite[Proposition 3.5]{DM2}, Lemma~\ref{Lemma 5.1} and Lemma~\ref{Lemma 5.2} that the holomorphic sectional
curvature $\tilde{H}_{t,n}$ of $D^{n}$ is given by
$$
\tilde{H}_{t,n}(E)=\kappa\|X\|^4+t\|A\|^4+\frac{(-1)^{n+1}}4(3+(-1)^{n+1}+4\kappa t)%
\|X\|^2\|A\|^2
$$
where $X=\pi _{*}E$, $A={\cal V}E$ and $\|E\|^2_{h_{t}}=\|X\|^2+t\|A\|^2 =1.$ Hence, for $n=1,$ the identity
$\tilde{H}_{t,n}\equiv\kappa$ is equivalent to $t=1/\kappa$, while for $n=2$ it is impossible. Thus the
proposition is proved.

\bth{Remark.} Similar arguments show that the Levi-Civita connection of the almost-Hermitian manifold $({\cal
Z},h_{t},J_{n})$ has a constant holomorphic sectional curvature $\kappa$ only in the case when $n=1$, $M$ is of
constant sectional curvature $\kappa $ and $t=1/\kappa  (\cite{DM2})$. \eth

\section{Examples of twistor spaces with parallel Nijenhuis tensor}

   It is well-known (\cite{M}) that the twistor space $({\cal Z},h_t,J_2)$ of an Einstein, self-dual manifold $M$
is a quasi-K\"ahler manifold satisfying the second Gray curvature condition. If $s>0$ and
$t=\displaystyle{\frac{6}{s}}$, resp. $s<0$ and $t=\displaystyle{-\frac{12}{s}}$ ($s$ is the scalar curvature of
$M$), then $({\cal Z},h_t,J_2)$ is nearly K\"ahler, resp. almost K\"ahler and, by results of \cite{BM} and
\cite{N}, the Nijenhuis tensor of $J_2$ is parallel with respect to the Chern connection. In fact, this is true
for any and $s$ and any $t$.

\begin{prop} Let $M$ be an Einstein and self-dual $4$-manifold with twistor space ${\cal Z}$. Then the Nijenhuis tensor of
the almost-complex structure $J_2$ is parallel with respect to the Chern connection of the almost-Hermitian
manifold $({\cal Z},h_{t},J_{2})$
\end{prop}

Proof. Denote by $N$ the Nijenhuis tensor of the almost-complex structure $J_2$. Let $\sigma$ be a point of ${\cal
Z}$ , $X,Y,Z$ vector fields on $M$ near the point $p=\pi(\sigma)$, and $A,B$ vertical vector fields near
${\sigma}$

The identity
$$
N(E,F)=-J_2(D_{E}J_2)(F)+J_2(D_{F}J_2)(E)-(D_{J_2F}J_2)(E)+(D_{J_2E}J_2)(F)
$$
and Lemma~\ref{Lemma 5.2} imply the following formulas:
\begin{equation}\label{eq 5.19}
 N(X^h,Y^h)_{\sigma}=-\frac s3(X\land K_{\sigma}Y+ K_{\sigma}X\land
Y);\> N(X^h,A)_{\sigma}=2(K_{\sigma\times A}X)_{\sigma}^h
\end{equation}

    As to the Chern connection $D^2$ of $({\cal Z},h_{t},J_{2})$, formulas (\ref{eq 5.1}), (\ref{eq 2.7}), (\ref{eq 2.8})
and Lemma~\ref{Lemma 5.2}  give:
\begin{equation}\label{eq 5.20}
\begin{array}{c}
D^{2}_{X^h}Y^h=(\nabla_{X}Y)^h;\>(D^{2}_{A}X^h)_{\sigma}= \displaystyle{\frac 12}(K_{\sigma\times A}X)^h;\\
                                                                                          \\
D^{2}_{X^h}A=D_{X^h}A-\displaystyle{\frac{ts}{24}}(K_{\sigma\times A}X)^h_{\sigma} ={\cal V}D_{X^h}A=[X^h,A].
\end{array}
\end{equation}

  Let $\xi$ be a section of ${\cal Z}$ near $p$
such that $\xi(p)=\sigma$ and $\nabla\xi|_p=0$. Then it is easy to see that, at the point $p$, $\nabla K_{\xi}=0$
and  ${\cal V}D_{X^h}A=\nabla_{X}(A\circ\xi)$ for any vertical vector field $A$ where $A\circ\xi$ is considered as
a section of $\Lambda^{2}_{-}TM$. Now the identities (\ref{eq 5.19}) and (\ref{eq 5.20}) imply that
$$
(D^{2}_{Z^h}N)(X^h,Y^h)_{\sigma}=0
$$

   Let $E_1,E_2,E_3,E_4$ be an oriented orthonormal frame of $TM$ near $p$
such that $\nabla E_i|_p=0, 1\leq i\leq 4$, and $s_{1}(p)=\sigma$ where $(s_1,s_2,s_3)$ is the local frame of
$\Lambda^{2}_{-}TM$ defined by (\ref{eq 2.1}). Then by (\ref{eq 5.19}) we have
$$
N(X^h,A)=-4\sum_{i=1}^4 g(J_2A,(X\wedge E_i)\circ\pi)E_i^h.
$$
Let us also note that
$$
Z^h_{\sigma}g(J_2A,(X\wedge E_i)\circ\pi)= Z_pg((J_2A)\circ\xi,X\wedge E_i)=
$$
$$
g(J_2\nabla_{Z_p}(A\circ\xi),X\wedge E_i)+g(J_2A_{\sigma}, \nabla_{Z_p}X\wedge E_i)
$$

    Now it is clear that
$$
(D^{2}_{Z^h}N)(X^h,A)_{\sigma}=0
$$

  We have also
$$
(D^{2}_{Z^h}N)(A,B)_{\sigma}=0
$$
in view of (\ref{eq 5.20}) and the fact that $N$ vanishes for any two vertical vectors.

  We shall further use the notations introduced at the end of the proof of Proposition~\ref{Proposition 5.3} .

   The fibres of ${\cal Z}$ are totally geodesic submanifolds, K\"ahlerian with respect to $J_2$, so the Chern
connection $D^{2}$ coincides with the Levi-Civita connection $D$ of $h_t$ for vertical vectors. Since $D_{U}U$ and
$D_{U}J_1U$ are vertical vectors and $[U,J_1U]_{\sigma}=0$, it follows from the standard formula for the
Levi-Civita connection that
\begin{equation}\label{eq 5.21}
(D_{U}U)_{\sigma}=(D_{U}J_1U)_{\sigma}=0
\end{equation}
 Hence
$$
D^{2}_{U_{\sigma}}N(X^h,Y^h)= U_{\sigma}(g(N(X^h,Y^h),U)U_{\sigma}+ U_{\sigma}(g(N(X^h,Y^h),J_1U))J_1U_{\sigma}.
$$
By (\ref{eq 5.20}) and (\ref{eq 2.6}), we have
$$
g(N(X^h,Y^h),U)=\frac{2s}{3}g(J_1U,X\wedge Y)=
$$
$$
\frac{2s}{3}(1-y_3)^{-1/2}(y_1y_3g(s_1,X\wedge Y)\circ\pi+ y_2y_3g(s_2,X\wedge Y)\circ\pi-(1-y_3^2)g(s_3,X\wedge
Y)\circ\pi)
$$
and
$$
g(N(X^h,Y^h),J_1U)=-\frac{2s}{3}g(U,X\wedge Y)=
$$
$$
\frac{2s}{3}(1-y_3)^{-1/2}(y_2g(s_1,X\wedge Y)\circ\pi- y_1g(s_2,X\wedge Y)\circ\pi)
$$
It follows that
$$
D_{U_{\sigma}}N(X^h,Y^h)=-\frac{2s}{3}g(s_1,X\wedge Y)_ps_3(p)
$$
Using (\ref{eq 5.20}), (\ref{eq 5.19}) and (\ref{eq 2.6C}) we easily obtain
$$
N(D^{2}_{U}X^h,Y^h)_{\sigma}=\frac{s}{6}g(X,Y)_ps_2(p)- \frac{s}{3}g(s_1,X\wedge Y)_ps_3(p),
$$
$$
N(X^h,D^{2}_{U}Y^h)_{\sigma}=-\frac{s}{6}g(X,Y)_ps_2(p)- \frac{s}{3}g(s_1,X\wedge Y)_ps_3(p).
$$
It follows that
$$
(D^{2}_{U}N)(X^h,Y^h)_{\sigma}=0
$$
Similarly
$$
(D^{2}_{J_1U}N)(X^h,Y^h)_{\sigma}=0.
$$
Therefore
$$
(D^{2}_{A}N)(X^h,Y^h)_{\sigma}=0
$$
for any vertical vector $A$ at $\sigma$.

  By (\ref{eq 5.19}) and (\ref{eq 5.20}) we obtain
$$
D^{2}_{U_{\sigma}}N(X^h,U)=-X^h_{\sigma} \>  \hbox   {and} \>
D^{2}_{U_{\sigma}}N(X^h,J_1U)=(K_{\sigma}X)^h_{\sigma}
$$
Taking also into account (\ref{eq 2.6C}), we get
$$
N(D^{2}_{U_{\sigma}}X^h,U)=-X^h_{\sigma} \>   \hbox   {and} \>
N(D^{2}_{U_{\sigma}}X^h,J_1U)=(K_{\sigma}X)^h_{\sigma}
$$
Then, by (\ref{eq 5.21}), we have
$$
(D^{2}_{U}N)(X^h,U)_{\sigma}=(D^{1}_{U}N)(X^h,J_1U)_{\sigma}=0
$$
Similarly, we get
$$
(D^{2}_{J_1U}N)(X^h,U)_{\sigma}=(D^{1}_{J_1U}N)(X^h,J_1U)_{\sigma}=0
$$
Therefore
$$
(D^{2}_{A}N)(X^h,B)_{\sigma}=0
$$

     Finally, let $A,B,C$ be vertical vectors at ${\sigma}$. Since
for vertical vectors $D^{1}$ coincides with the Levi-Civita connection of the fibre through $\sigma$, we have
$$
(D^{2}_{A}N)(B,C)_{\sigma}=0
$$

{\bf Remark.}  Identity (\ref{eq 5.20}) shows that the Chern connection $D^2$ actually does not depend on $t$ when
the base manifold is Einstein and self-dual. Then Proposition 4 for $s\neq 0$ follows also from the results in
\cite{BM} and \cite{N} mentioned in the beginning of this section. Propositions 1 and 4 can be extended to the
twistor spaces of quaternionic-K\"ahler manifolds by means of the formulas in \cite{AGI}.

\end{document}